\documentclass{amsart}
\usepackage{amssymb}
\usepackage{amsmath}
\usepackage{amsfonts}
\usepackage{graphicx,color}

\newtheorem{theorem}{Theorem}
\theoremstyle{plain}

\newtheorem{lemma}{Lemma}

\numberwithin{equation}{section}
\input{tcilatex}

\begin{document}
\title[Singular quasilinear elliptic systems]{Singular quasilinear elliptic
systems\\
with (super-) homogeneous condition}
\author{Hana Didi$^{\ast },$ Brahim Khodja $^{\ast }$}
\address{$\ast ${\small \textit{Mathematic Department, Badji-Mokhtar Annaba
University, }}\\
\textit{23}{\small \textit{000 Annaba Algeria}}}
\email{hana.di@hotmail.fr, \ brahim.khodja@univ-annaba.org}
\author{Abdelkrim Moussaoui$^{\ast \ast }$}
\address{$\ast \ast $ {\small \textit{Biology Department, A. Mira Bejaia
University, }}\\
{\small \textit{Targa Ouzemour 06000 Bejaia, Algeria}}}
\email{abdelkrim.moussaoui@univ-bejaia.dz}
\thanks{}

\begin{abstract}
In this paper we establish existence, nonexistence and regularity of
positive solutions for a class of singular quasilinear elliptic systems
subject to (super-) homogeneous condition. The approach is based on
sub-supersolution methods for systems of quasilinear singular equations
combined with perturbation arguments involving singular terms.
\end{abstract}

\maketitle

\section{Introduction and main results\label{S1}}

We consider the following system of quasilinear and singular elliptic
equations:%
\begin{equation*}
(\mathcal{P})\qquad \left\{ 
\begin{array}{l}
-\Delta _{p_{1}}u_{1}=\lambda u_{1}^{\alpha _{1}}u_{2}^{\beta _{1}}+\delta
h_{1}(x)\text{ in }\Omega \\ 
-\Delta _{p_{2}}u_{2}=\lambda u_{1}^{\alpha _{2}}u_{2}^{\beta _{2}}+\delta
h_{2}(x)\text{ in }\Omega \\ 
u_{1},u_{2}>0\text{ in }\Omega \\ 
u_{1},u_{2}=0\text{ on }\partial \Omega ,%
\end{array}%
\right.
\end{equation*}%
where $\Omega $ is a bounded domain in $%
\mathbb{R}
^{N}$ $\left( N\geq 2\right) $ having a smooth boundary $\partial \Omega $, $%
\lambda >0,$ $\delta \geq 0$ are parameters, and $h_{i}\in L^{\infty
}(\Omega )$ is a nonnegative function. Here $\Delta _{p_{i}}$ ($1<p_{i}\leq
N $) stands for the $p_{i}$-Laplacian operator. We consider the system ($%
\mathcal{P}$) in a singular case assuming that%
\begin{equation}
\begin{array}{c}
0<\alpha _{2}<p_{1}^{\ast }-1\text{, \ }0<\beta _{1}<p_{2}^{\ast }-1\text{
and }-2<\alpha _{1},\beta _{2}<0\text{,}%
\end{array}
\label{h1}
\end{equation}%
where $p_{i}^{\ast }=\frac{Np_{i}}{N-p_{i}}$. This assumption makes system ($%
\mathcal{P}$) be cooperative, that is, for $u_{1}$ (resp. $u_{2}$) fixed the
right term in the first (resp. second) equation of ($\mathcal{P}$) is
increasing in $u_{2}$ (resp. $u_{1}$).

Recently, singular cooperative system ($\mathcal{P}$) with $\delta =0$ was
mainly studied in \cite{GHM,GHS,MM}. In \cite{MM} existence and boundedness
theorems for ($\mathcal{P}$) was established by using sub-supersolution
method for systems combined with perturbation techniques. In \cite{GHM} one
gets existence, uniqueness, and regularity of a positive solution on the
basis of an iterative scheme constructed through a sub-supersolution pair.
In \cite{GHS} an existence theorem involving sub-supersolution was obtained
through a fixed point argument in a sub-supersolution setting. The
semilinear case in ($\mathcal{P}$) (i.e. $p_{i}=2$) was considered in \cite%
{G,HMV,MKT} where the linearity of the principal part is essentially used.
In this context, the singular system ($\mathcal{P}$) can be viewed as the
elliptic counter-part of a class of Gierer-Meinhardt systems that models
some biochemical processes (see, e.g. \cite{MKT}). It can be also given an
astrophysical meaning since it generalize to system the well-known
Lane-Emden equation, where all exponents are negative (see \cite{G}). The
complementary situation for system ($\mathcal{P}$) with respect to (\ref{h1}%
) is the so-called competitive system, which has recently attracted much
interest. Relevant contributions regarding this topic can be found in \cite%
{GHS,MM2,MM3}. It is worth pointing out that the aforementioned works have
examined the subhomogeneous case $\Theta >0$ of singular problem ($\mathcal{P%
}$) where%
\begin{equation}
\begin{array}{c}
\Theta =\left( p_{1}-1-\alpha _{1}\right) \left( p_{2}-1-\beta _{2}\right)
-\beta _{1}\alpha _{2}.%
\end{array}
\label{c0}
\end{equation}%
The constant $\Theta $ is related to system stability ($\mathcal{P}$) that
behaves in a drastically different way, depending on the sign of $\Theta $.
For instance, for $\Theta <0$ system ($\mathcal{P}$) is not stable in the
sense that possible solutions cannot be obtained by iterative methods (see 
\cite{CFMT}).

Unlike the subhomogeneous case $\Theta >0$ studied in the above references,
the novelty of this paper is to establish the existence, regularity and
nonexistence of (positive) solutions for singular problem ($\mathcal{P}$) by
processing the two cases: 'homogeneous' when $\Theta =0$ and
'superhomogeneous' if $\Theta <0$. It should be noted that throughout this
paper, $\Theta <0$ (resp. $=0$) means that $p_{i}-1-\alpha _{i}-\beta _{i}<0$
(resp. $=0$).

The existence result for problem ($\mathcal{P}$) is stated as follows.

\begin{theorem}
\label{T4}Assume {(\ref{h1})},{\ }$\Theta <0$ (resp. $\Theta =0$) and
suppose that%
\begin{equation}
\inf_{\Omega }h_{1}(x),\text{ }\inf_{\Omega }h_{2}(x)>0.  \label{c1}
\end{equation}%
Then, there is $\delta _{0}>0$ (resp. $\delta _{0},\lambda _{0}>0$) such
that, for all $\delta \in (0,\delta _{0})$, problem ($\mathcal{P}$)
possesses a (positive) solution $(u_{1},u_{2})$ in $(W_{0}^{1,p_{1}}(\Omega
)\cap L^{\infty }(\Omega ))\times (W_{0}^{1,p_{2}}(\Omega )\cap L^{\infty
}(\Omega ))$ verifying $u_{i}\geq cd(x)$ in $\Omega $, for some constant $%
c>0 $ and for all $\lambda >0$ (resp. $\lambda \in (0,\lambda _{0})$).
Moreover, if $\alpha _{1},\beta _{2}>-1$, the solution $(u_{1},u_{2})$ is
bounded in $C_{0}^{1,\beta }(\overline{\Omega })\times C_{0}^{1,\beta }(%
\overline{\Omega }),$ for certain $\beta \in (0,1)$. Furthermore, if $\Theta
=\delta =0$ and 
\begin{equation}
\begin{array}{l}
\beta _{1}=\frac{p_{2}}{p_{1}}(p_{1}-1-\alpha _{1})\text{ \ or \ }\alpha
_{2}=\frac{p_{1}}{p_{2}}(p_{2}-1-\beta _{2}),%
\end{array}
\label{c2}
\end{equation}%
then, there exists $\lambda _{\ast }>0$ such that problem ($\mathcal{P}$)
has no solution for every $\lambda \in (0,\lambda _{\ast })$.
\end{theorem}

\bigskip

A solution of ($\mathcal{P}$) is understood in the weak sense, that is, a
pair $(u_{1},u_{2})\in W_{0}^{1,p_{1}}(\Omega )\times W_{0}^{1,p_{2}}(\Omega
)$, which are positive a.e. in $\Omega $ and satisfying 
\begin{equation*}
\int_{\Omega }|\nabla u_{i}|^{p_{i}-2}\nabla u_{i}\nabla \varphi _{i}\
dx=\int_{\Omega }(\lambda u_{1}^{\alpha _{i}}u_{2}^{\beta _{i}}+\delta
h_{i})\varphi _{i}\ dx,\text{\ for all }\varphi _{i}\in
W_{0}^{1,p_{i}}(\Omega ),\text{ }i=1,2.
\end{equation*}

\bigskip

The main technical{\small \ }difficulty consists in the presence{\small \ }%
of singular terms in system ($\mathcal{P}$) with (\ref{h1}), expressed
through (super-) homogeneous condition. Our approach is chiefly based on
sub-supersolution method in its version for systems \cite[section 5.5]{CLM}.
However, this method cannot be directly implemented due to the presence of
singular terms in ($\mathcal{P}$) under assumption (\ref{h1}). So, we first
disturb system ($\mathcal{P}$) by introducing a parameter $\varepsilon >0$.
This gives rise to a regularized system for ($\mathcal{P}$) depending on $%
\varepsilon $ whose study is relevant for our initial problem. By applying
the sub-supersolution method, we show that the regularized system has a
positive solution $(u_{1,\varepsilon },u_{2,\varepsilon })$ in $C^{1,\beta }(%
\overline{\Omega })\times C^{1,\beta }(\overline{\Omega })$ for some $\beta
\in (0,1)$. It is worth noting that the choice of suitable functions with an
adjustment of adequate constants is crucial in order to construct the
sub-supersolution pair as well as to process the both cases $\Theta <0$ and $%
\Theta =0$. The (positive) solution $(u_{1},u_{2})$ in $(W_{0}^{1,p_{1}}(%
\Omega )\cap L^{\infty }(\Omega ))\times (W_{0}^{1,p_{2}}(\Omega )\cap
L^{\infty }(\Omega ))$ of ($\mathcal{P}$) is obtained by passing to the
limit as $\varepsilon \rightarrow 0$. This is based on a priori estimates,
Fatou's Lemma and $S_{+}$-property of the negative $p_{i}$-Laplacian. The
positivity of the solution $(u_{1},u_{2})$ is achieved through assumption (%
\ref{c1}) while $C^{1,\beta }$-regularity is derived from the regularity
result in \cite{Hai}.

The rest of the paper is organized as follows. Section \ref{S2} is devoted
to the existence of solutions for the regularized system. Section \ref{S3}
established the proof of the main result.

\section{The regularized system\label{S2}}

Given $1<p<+\infty $, the space $L^{p}(\Omega )$ and $W_{0}^{1,p}(\Omega )$
are endowed with the usual norms $\left\Vert u\right\Vert _{p}=(\int_{\Omega
}\left\vert u\right\vert ^{p}dx)^{1/p}$ and $\left\Vert u\right\Vert
_{1,p}=(\int_{\Omega }\left\vert \nabla u\right\vert ^{p}dx)^{1/p},$
respectively. We will also utilize the space $C_{0}^{1,\beta }(\overline{%
\Omega })=\left\{ u\in C^{1,\beta }(\overline{\Omega }):u=0\text{ on }%
\partial \Omega \right\} $ for a suitable $\beta \in (0,1)$. 

In what follows, we denote by $\phi _{1,p_{i}}$ the positive eigenfunction
associated with the principal eigenvalue $\lambda _{1,p_{i}},$ characterized
by the minimum of Rayleigh quotient%
\begin{equation}
\lambda _{1,p_{i}}=\inf_{u_{i}\in W_{0}^{1,p_{i}}\left( \Omega \right)
\backslash \{0\}}\frac{\left\Vert \nabla u_{i}\right\Vert _{p_{i}}^{p_{i}}}{%
\left\Vert u_{i}\right\Vert _{p_{i}}^{p_{i}}}.  \label{62}
\end{equation}%
For a later use recall there exist constants $l_{i,}\hat{l}_{i}>0$ such that%
\begin{equation}
\hat{l}_{1}\phi _{1,p_{1}}(x)\geq \phi _{1,p_{2}}(x)\geq \hat{l}_{2}\phi
_{1,p_{1}}(x)\text{ and }l_{1}d(x)\geq \phi _{1,p_{i}}(x)\geq l_{2}d(x)\text{
for all }x\in \Omega   \label{6}
\end{equation}%
where $d(x):=dist(x,\partial \Omega )$ (see, e.g., \cite{GST}).

Let $\tilde{\Omega}$ be a bounded domain in $%
\mathbb{R}
^{N}$ with a smooth boundary $\partial \tilde{\Omega}$ such that $\overline{%
\Omega }\subset \tilde{\Omega}.$ Denote $\tilde{d}(x):=d(x,\partial \tilde{%
\Omega})$. By the definition of $\tilde{\Omega}$ there exists a constant $%
\rho >0$ sufficiently small such that%
\begin{equation}
\tilde{d}(x)>\rho \text{ in }\overline{\Omega }.  \label{10}
\end{equation}%
Define $w_{i}\in C^{1}(\overline{\tilde{\Omega}})$ the unique solution of
the torsion problem 
\begin{equation}
-\Delta _{p_{i}}w_{i}=1\text{ in }\tilde{\Omega},\text{ }w_{i}=0\text{ on }%
\partial \tilde{\Omega},  \label{12}
\end{equation}%
satisfying the estimates 
\begin{equation}
w_{i}(x)\geq c_{0}\tilde{d}(x)\text{\ in}\ \tilde{\Omega}\text{, for certain 
}c_{0}\in (0,1).  \label{11}
\end{equation}%
For a real constant $C>1$, set 
\begin{equation}
(\underline{u}_{i,\varepsilon },\overline{u}_{i})=(c_{\varepsilon }\phi
_{1,p_{i}},C^{-1}w_{i}),\text{ }i=1,2,  \label{15}
\end{equation}%
where $c_{\varepsilon }>0$ is a constant depending on $\varepsilon >0$ such
that 
\begin{equation}
0<c_{\varepsilon }<c_{0}l_{1}^{-1}C^{-1}.  \label{13}
\end{equation}%
Then, by (\ref{15}), (\ref{6}) and (\ref{12}), it is readily seen that%
\begin{equation*}
\begin{array}{l}
\overline{u}_{i}(x)=C^{-1}w_{i}(x)\geq C^{-1}c_{0}\tilde{d}(x)\geq
C^{-1}c_{0}d(x) \\ 
\geq l_{1}^{-1}C^{-1}c_{0}\phi _{1,p_{i}}(x)\geq c_{\varepsilon }\phi
_{1,p_{i}}(x)=\underline{u}_{i,\varepsilon }(x)\text{ in }\overline{\Omega },%
\text{ for }i=1,2.%
\end{array}%
\end{equation*}

For every $\varepsilon \in (0,\varepsilon _{0})$, with $\varepsilon _{0}<1,$
let introduce the auxiliary problem%
\begin{equation*}
(\mathcal{P}_{\varepsilon })\qquad \left\{ 
\begin{array}{ll}
-\Delta _{p_{1}}u_{1}=\lambda (u_{1}+\varepsilon )^{\alpha
_{1}}(u_{2}+\varepsilon )^{\beta _{1}}+\delta h_{1}(x) & \text{in }\Omega \\ 
-\Delta _{p2}u_{2}=\lambda (u_{1}+\varepsilon )^{\alpha
_{2}}(u_{2}+\varepsilon )^{\beta _{2}}+\delta h_{2}(x) & \text{in }\Omega \\ 
u_{1},u_{2}=0 & \text{on }\partial \Omega ,%
\end{array}%
\right.
\end{equation*}%
which provides approximate solutions for the initial problem ($\mathcal{P}$).

\begin{lemma}
\label{L2}Assume (\ref{h1}) and $h_{1},h_{2}\neq 0$ in $\Omega $. Then, if $%
\Theta <0$ (resp. $\Theta =0$), there is a constant $\delta _{0}>0$ (resp. $%
\delta _{0},\lambda _{0}>0$ ) such that for all $\delta \in (0,\delta _{0}),$
$(\overline{u}_{1},\overline{u}_{2})$ in (\ref{15}) is a supersolution of ($%
\mathcal{P}_{\varepsilon }$) for all $\lambda >0$ (resp. $\lambda \in
(0,\lambda _{0})$) and all $\varepsilon \in (0,\varepsilon _{0})$.
\end{lemma}

\begin{proof}
Assume $\Theta <0$ and set $\varepsilon _{0}=C^{-1},$%
\begin{equation}
\delta _{0}=\frac{1}{2}\min_{i=1,2}\{\frac{1}{C^{p_{i}-1}\left\Vert
h_{i}\right\Vert _{\infty }}\}.  \label{33}
\end{equation}%
On account of (\ref{h1}), (\ref{10})-(\ref{15}) and (\ref{33}), for all $%
\delta \in (0,\delta _{0})$ and $\varepsilon \in (0,\varepsilon _{0}),$ one
derives%
\begin{equation*}
\begin{array}{l}
(\overline{u}_{1}+\varepsilon )^{-\alpha _{1}}(\overline{u}_{2}+\varepsilon
)^{-\beta _{1}}(-\Delta _{p_{1}}\overline{u}_{1}-\delta h_{1})\geq \overline{%
u}_{1}^{-\alpha _{1}}(\overline{u}_{2}+\varepsilon _{0})^{-\beta
_{1}}(-\Delta _{p_{1}}\overline{u}_{1}-\delta \left\Vert h_{1}\right\Vert
_{\infty }) \\ 
\geq C^{\alpha _{1}+\beta _{1}}(c_{0}\tilde{d}(x))^{-\alpha _{1}}(\left\Vert
w_{2}\right\Vert _{\infty }+1)^{-\beta _{1}}(C^{-(p_{1}-1)}-\delta
_{0}\left\Vert h_{1}\right\Vert _{\infty }) \\ 
\geq C^{\beta _{1}-(p_{1}-1-\alpha _{1})}(c_{0}\rho )^{-\alpha
_{1}}(\left\Vert w_{2}\right\Vert _{\infty }+1)^{-\beta _{1}}(1-\delta
_{0}C^{p_{1}-1}\left\Vert h_{1}\right\Vert _{\infty }) \\ 
\geq \frac{1}{2}C^{\beta _{1}-(p_{1}-1-\alpha _{1})}(c_{0}\rho )^{-\alpha
_{1}}(\left\Vert w_{2}\right\Vert _{\infty }+1)^{-\beta _{1}}\geq \lambda 
\text{ \ in }\overline{\Omega },%
\end{array}%
\end{equation*}%
and similarly 
\begin{equation*}
\begin{array}{l}
(\overline{u}_{1}+\varepsilon )^{-\alpha _{2}}(\overline{u}_{2}+\varepsilon
)^{-\beta _{2}}(-\Delta _{p_{2}}\overline{u}_{2}-\delta h_{2})\geq (%
\overline{u}_{1}+\varepsilon _{0})^{-\alpha _{2}}\overline{u}_{2}^{-\beta
_{2}}(-\Delta _{p_{2}}\overline{u}_{2}-\delta \left\Vert h_{2}\right\Vert
_{\infty }) \\ 
\geq C^{\alpha _{2}+\beta _{2}}(\left\Vert w_{1}\right\Vert _{\infty
}+1)^{-\alpha _{2}}(c_{0}\tilde{d}(x))^{-\beta _{2}}(C^{-(p_{2}-1)}-\delta
_{0}\left\Vert h_{2}\right\Vert _{\infty }) \\ 
=C^{\alpha _{2}-(p_{2}-1-\beta _{2})}(\left\Vert w_{1}\right\Vert _{\infty
}+1)^{-\alpha _{2}}(c_{0}\rho )^{-\beta _{2}}(1-\delta
_{0}C^{p_{2}-1}\left\Vert h_{2}\right\Vert _{\infty }) \\ 
\geq \frac{1}{2}C^{\alpha _{2}-(p_{2}-1-\beta _{2})}(\left\Vert
w_{1}\right\Vert _{\infty }+1)^{-\alpha _{2}}(c_{0}\rho )^{-\beta _{2}}\geq
\lambda \text{ \ in }\overline{\Omega },%
\end{array}%
\end{equation*}%
for all $\lambda >0$, provided $C>1$ is sufficiently large. This shows that $%
(\overline{u}_{1},\overline{u}_{2})$ is a supersolution pair for problem ($%
\mathcal{P}_{\varepsilon }$). If $\Theta =0$, by repeating the argument
above, the same conclusion can be drawn for $\lambda \in (0,\lambda _{0})$
with a constant $\lambda _{0}>0$ that can be precisely estimated. This
completes the proof.
\end{proof}

\begin{lemma}
\label{L1}Assume (\ref{h1}) and $\Theta \leq 0$ hold. Then, $(\underline{u}%
_{1,\varepsilon },\underline{u}_{2,\varepsilon })$ is a subsolution of ($%
\mathcal{P}_{\varepsilon }$) for all $\lambda ,\delta >0$ and every $%
\varepsilon \in (0,\varepsilon _{0})$.
\end{lemma}

\begin{proof}
Fix $\varepsilon \in (0,\varepsilon _{0})$. From (\ref{15}) and (\ref{h1}),
we obtain 
\begin{equation}
\begin{array}{l}
(\underline{u}_{1,\varepsilon }+\varepsilon )^{-\alpha _{1}}(\underline{u}%
_{2,\varepsilon }+\varepsilon )^{-\beta _{1}}(-\Delta _{p_{1}}\underline{u}%
_{1,\varepsilon }-\delta h_{1}) \\ 
\leq c_{\varepsilon }^{p_{1}-1}(c_{\varepsilon }\phi _{1,p_{1}}+\varepsilon
_{0})^{-\alpha _{1}}(c_{\varepsilon }\phi _{1,p_{2}}+\varepsilon )^{-\beta
_{1}}\lambda _{1,p_{1}}\phi _{1,p_{1}}^{p_{1}-1} \\ 
\leq c_{\varepsilon }^{p_{1}-1}(\phi _{1,p_{1}}+\varepsilon _{0})^{-\alpha
_{1}}(c_{\varepsilon }\phi _{1,p_{2}}+\varepsilon )^{-\beta _{1}}\lambda
_{1,p_{1}}\phi _{1,p_{1}}^{p-1} \\ 
\leq c_{\varepsilon }^{p_{1}-1}\varepsilon ^{-\beta _{1}}(\left\Vert \phi
_{1,p_{1}}\right\Vert _{\infty }+1)^{-\alpha _{1}}\lambda
_{1,p_{1}}\left\Vert \phi _{1,p_{1}}\right\Vert _{\infty }^{p_{1}-1}\leq
\lambda \text{ \ in }\overline{\Omega }%
\end{array}
\label{20}
\end{equation}%
and similarly%
\begin{equation}
\begin{array}{l}
(\underline{u}_{1,\varepsilon }+\varepsilon )^{-\alpha _{2}}(\underline{u}%
_{2,\varepsilon }+\varepsilon )^{-\beta _{2}}(-\Delta _{p_{2}}\underline{u}%
_{2,\varepsilon }-\delta h_{2}) \\ 
\leq (\underline{u}_{1,\varepsilon }+\varepsilon )^{-\alpha _{2}}(\underline{%
u}_{2,\varepsilon }+\varepsilon _{0})^{-\beta _{2}}(\Delta _{p_{2}}%
\underline{u}_{2,\varepsilon }) \\ 
=c_{\varepsilon }^{p_{2}-1}(c_{\varepsilon }\phi _{1,p_{1}}+\varepsilon
)^{-\alpha _{2}}(\phi _{1,p_{2}}+\varepsilon _{0})^{-\beta _{2}}\lambda
_{1,p_{2}}\phi _{1,p_{2}}^{p_{2}-1} \\ 
\leq c_{\varepsilon }^{p_{2}-1}\varepsilon ^{-\alpha _{2}}(\left\Vert \phi
_{1,p_{2}}\right\Vert _{\infty }+\varepsilon _{0})^{-\beta _{2}}\lambda
_{1,p_{2}}\left\Vert \phi _{1,p_{2}}\right\Vert _{\infty }^{p_{2}-1}\leq
\lambda \text{ \ in }\overline{\Omega },%
\end{array}
\label{21}
\end{equation}%
provided $c_{\varepsilon }>0$ is sufficiently small. Gathering (\ref{20})
and (\ref{21}) together yields 
\begin{equation*}
-\Delta _{p_{i}}\underline{u}_{i,\varepsilon }\leq \lambda (\underline{u}%
_{1,\varepsilon }+\varepsilon )^{\alpha _{i}}(\underline{u}_{2,\varepsilon
}+\varepsilon )^{\beta _{i}}+\delta h_{i}\ \ \text{in }\overline{\Omega },
\end{equation*}%
proving that $(\underline{u}_{1,\varepsilon },\underline{u}_{2,\varepsilon
}) $ in (\ref{15}) is a subsolution pair for problem ($\mathcal{P}%
_{\varepsilon }$).
\end{proof}

\bigskip

We state the following result regarding the regularized system.

\begin{theorem}
\label{T5}Assume (\ref{h1}) and $h_{1},h_{2}\neq 0$ in $\Omega $. Then
\end{theorem}

\begin{description}
\item[\textrm{(a)}] \textit{If }$\Theta <0$\textit{\ (resp. }$\Theta =0$) 
\textit{there exist a constant }$\delta _{0}>0$\textit{\ (resp. }$\delta
_{0},\lambda _{0}>0$)\textit{\ such that for all }$\delta \in (0,\delta
_{0}) $ \textit{system (}$\mathcal{P}_{\varepsilon }$\textit{) has a
(positive) solution }$(u_{1,\varepsilon },u_{2,\varepsilon })\in
C_{0}^{1,\beta }(\overline{\Omega })\times C_{0}^{1,\beta }(\overline{\Omega 
}),$\textit{\ }$\beta \in (0,1)$\textit{, satisfying }%
\begin{equation}
u_{i,\varepsilon }(x)\leq \overline{u}_{i}(x)\text{ in }\Omega ,  \label{22}
\end{equation}%
for all$\mathit{\ }\lambda >0$ (resp. $\lambda \in (0,\lambda _{0})$), and
every$\mathit{\ }\varepsilon \in (0,\varepsilon _{0})$.

\item[\textrm{(b)}] \textit{For }$\Theta \leq 0$\textit{\ and under
assumption (\ref{c1}), if }$\delta >0$, \textit{there exists a constant }$%
c_{0}>0$\textit{, independent of }$\varepsilon $\textit{, such that all
solutions }$(u_{1,\varepsilon },u_{2,\varepsilon })$\textit{\ of system (}$%
\mathcal{P}_{\varepsilon }$\textit{) verify }%
\begin{equation}
u_{i,\varepsilon }(x)\geq c_{0}d(x)\text{ for a.a. }x\in \Omega ,\text{ for
all }\varepsilon \in (0,\varepsilon _{0}).  \label{23}
\end{equation}
\end{description}

\begin{proof}
On the basis of Lemmas \ref{L2} and \ref{L1} together with \cite[section 5.5]%
{CLM} there exists a solution $\left( u_{1,\varepsilon },u_{2,\varepsilon
}\right) $ of problem ($\mathcal{P}_{\varepsilon }$), for every $\varepsilon
\in (0,\varepsilon _{0})$. Moreover, applying the regularity theory (see 
\cite{L}), we infer that $\left( u_{1,\varepsilon },u_{2,\varepsilon
}\right) \in C_{0}^{1,\beta }(\overline{\Omega })\times C_{0}^{1,\beta }(%
\overline{\Omega })$ for a suitable $\beta \in (0,1)$. This proves \textrm{%
(a)}.

Now, according to (\ref{c1}), let $\sigma >0$ be a constant such that $%
\inf_{\Omega }h_{1}(x),$ $\inf_{\Omega }h_{2}(x)>\sigma $. Define $z_{i}$
the only positive solution of 
\begin{equation*}
\begin{array}{l}
-\Delta _{p_{i}}z_{i}=\delta \sigma \text{ \ in }\Omega ,\text{ }z_{i}=0%
\text{ \ on }\partial \Omega ,%
\end{array}%
\end{equation*}
which is known to satisfy $z_{i}(x)\geq c_{2}d(x)$ in $\Omega $. Then it
follows that $-\Delta _{p_{i}}u_{\varepsilon }\geq -\Delta _{p_{i}}z_{i}$ in 
$\Omega ,$ \ $u_{i,\varepsilon }=z_{i}$ on $\partial \Omega ,$\ for all $%
\varepsilon \in (0,\varepsilon _{0}),$ and therefore, the weak comparison
principle ensures the assertion \textrm{(b) }holds true.
\end{proof}

\section{Proof of the main result\label{S3}}

Set $\varepsilon =\frac{1}{n}$ with any positive integer $n>1/\varepsilon
_{0}$. From Theorem \ref{T5} with $\varepsilon =\frac{1}{n}$, there exists $%
u_{i,n}:=u_{i,\frac{1}{n}}$ such that%
\begin{equation}
\begin{array}{l}
\langle -\Delta _{p_{i}}u_{i,n},\varphi _{i}\rangle =\lambda \int_{\Omega
}(u_{1,n}+\frac{1}{n})^{\alpha _{i}}(u_{2,n}+\frac{1}{n})^{\beta
_{i}}\varphi _{i}\ dx+\delta \int_{\Omega }h_{i}\varphi _{i}dx,%
\end{array}
\label{122}
\end{equation}%
for all $\varphi _{i}\in W_{0}^{1,p_{i}}(\Omega ),i=1,2$. Taking $\varphi
=u_{1,n}$ in (\ref{122}), since $\alpha _{1}<0<\beta _{1}$, we get 
\begin{equation}
\Vert u_{1,n}\Vert _{1,p_{1}}^{p_{1}}\leq \lambda \int_{\Omega
}u_{1,n}^{\alpha _{1}+1}\left( u_{2,n}+1\right) ^{\beta _{1}}dx+\delta
\left\Vert h_{1}\right\Vert _{\infty }\int_{\Omega }u_{1,n}dx.  \label{123}
\end{equation}%
If $-1\leq \alpha _{1}<0$ (see (\ref{h1})), on the basis of (\ref{22}) with $%
\varepsilon =\frac{1}{n}$, it follows directly that $\{u_{1,n}\}$ is bounded
in $W_{0}^{1,p_{1}}\left( \Omega \right) $. It remains to argue when $%
-2<\alpha _{1}<-1$. Here, we point out that the subsolution $(\underline{u}%
_{1,\varepsilon },\underline{u}_{2,\varepsilon })$ depends on $\varepsilon $
(see Lemma \ref{L1}) and therefore it is inoperable. However, this
difficulty is avoided by taking advantage of (\ref{22}) and (\ref{23}),
which lead to 
\begin{equation*}
\begin{array}{c}
\Vert u_{1,n}\Vert _{1,p_{1}}^{p_{1}}\leq C_{0}\int_{\Omega }d(x)^{\alpha
_{1}+1}\ dx+C_{0}^{\prime }\delta \left\Vert h_{1}\right\Vert _{\infty
}\left\Vert u_{1,n}\right\Vert _{1,p_{1}},%
\end{array}%
\end{equation*}%
for certain constants $C_{0}$ and $C_{0}^{\prime }$ independent of $n$.
Since $\alpha _{1}+1>-1$ and $p_{1}>1$, we infer from \cite[Lemma in page 726%
]{LM} that $\{u_{1,n}\}$ is bounded in $W_{0}^{1,p_{1}}\left( \Omega \right) 
$. Similarly, we derive that $\{u_{2,n}\}$ is bounded in $%
W_{0}^{1,p_{2}}\left( \Omega \right) $. We are thus allowed to extract
subsequences (still denoted by $\{u_{i,n}\}$) such that 
\begin{equation}
\begin{array}{c}
u_{i,n}\rightharpoonup u_{i}\text{ in }W_{0}^{1,p_{i}}(\Omega )\text{, }%
i=1,2.%
\end{array}
\label{30}
\end{equation}%
The convergence in (\ref{30}) combined with Rellich embedding Theorem and (%
\ref{22})-(\ref{23}) entails 
\begin{equation}
c_{0}d(x)\leq u_{i}(x)\leq \overline{u}_{i}(x)\ \ \text{in }\Omega .
\label{320}
\end{equation}%
Inserting $\varphi _{i}=u_{i,n}-u_{i}$ in (\ref{122}) yields \ 
\begin{equation*}
\langle -\Delta _{p_{i}}u_{i,n},u_{i,n}-u_{i}\rangle =\int_{\Omega }[\lambda
(u_{1,n}+\frac{1}{n})^{\alpha _{i}}(u_{2,n}+\frac{1}{n})^{\beta _{i}}+\delta
h_{i}](u_{i,n}-u_{i})\ dx.
\end{equation*}%
We claim that $\underset{n\rightarrow \infty }{\lim }\langle -\Delta
_{p_{i}}u_{i,n},u_{i,n}-u_{i}\rangle \leq 0.$ Indeed, from (\ref{22}), (\ref%
{23}) and (\ref{15}), we have%
\begin{equation*}
\begin{array}{c}
\lbrack \left( u_{1,n}+\frac{1}{n}\right) ^{\alpha _{1}}\left( u_{2,n}+\frac{%
1}{n}\right) ^{\beta _{1}}+\delta h_{1}](u_{1,n}-u_{1})\leq \left[
u_{1,n}^{\alpha _{1}}\left( u_{2,n}+1\right) ^{\beta _{1}}+\delta h_{1}%
\right] (u_{1,n}-u_{1}) \\ 
\leq u_{1,n}^{\alpha _{1}+1}\left( u_{2,n}+1\right) ^{\beta _{1}}+\delta
h_{1}u_{1,n}\leq \Phi (x),%
\end{array}%
\end{equation*}%
where \cite[Lemma in page 726]{LM}, with $\alpha _{1}>-2,$ garantees the
integrability of function $\Phi $ defined by%
\begin{equation*}
\Phi (x)=\left\{ 
\begin{array}{ll}
(C^{-1}\left\Vert w_{1}\right\Vert _{\infty })^{\alpha
_{1}+1}(C^{-1}\left\Vert w_{2}\right\Vert _{\infty }+1)^{\beta _{1}}+\delta
C^{-1}\left\Vert h_{1}\right\Vert _{\infty }\left\Vert w_{1}\right\Vert
_{\infty } & \text{if }\alpha _{1}\geq -1 \\ 
(c_{2}d(x))^{\alpha _{1}+1}(C^{-1}\left\Vert w_{2}\right\Vert _{\infty
}+1)^{\beta _{1}}+\delta C^{-1}\left\Vert h_{1}\right\Vert _{\infty
}\left\Vert w_{1}\right\Vert _{\infty } & \text{if }\alpha _{1}<-1.%
\end{array}%
\right.
\end{equation*}%
By (\ref{30}) and applying Fatou's Lemma, the claim follows. A similiar
argument shows that $\underset{n\rightarrow \infty }{\lim }\left\langle
-\Delta _{p_{2}}u_{2,n},u_{2,n}-u_{2}\right\rangle \leq 0$. Then the $S_{+}$%
-property of $-\Delta _{p_{i}}$ on $W_{0}^{1,p_{i}}(\Omega )$ guarantees that%
\begin{equation*}
u_{i,n}\longrightarrow u_{i}\text{ \ in }W_{0}^{1,p_{i}}(\Omega ),\text{ }%
i=1,2.
\end{equation*}
Hence we may pass to the limit in (\ref{122}) to conclude that $%
(u_{1},u_{2}) $ is a solution of problem ($\mathcal{P}$) satsifying (\ref%
{320}). Furthermore, using (\ref{h1}), (\ref{320}) and (\ref{15}), one has%
\begin{equation*}
u_{1}^{\alpha _{1}}u_{2}^{\beta _{1}}+\delta h_{1}\leq C_{1}^{\prime
}d(x)^{\alpha _{1}}\text{ and }u_{1}^{\alpha _{2}}u_{2}^{\beta _{2}}+\delta
h_{2}\leq C_{2}^{\prime }d(x)^{\beta _{2}}\text{ in }\Omega ,
\end{equation*}%
for certain constants $C_{1}^{\prime },C_{2}^{\prime }>0$. Then, for $\alpha
_{1},\beta _{2}>-1$ and owing to \cite[Lemma 3.1]{Hai}, we infer that $%
(u_{1},u_{2})\in C_{0}^{1,\beta }(\overline{\Omega })\times C_{0}^{1,\beta }(%
\overline{\Omega })$ for certain $\beta \in (0,1)$.

We are left with the task of determining the nonexistence result stated in
Theorem \ref{T4}. Arguing by contradiction and assume that $(u_{1},u_{2})$
is a positive solution of problem ($\mathcal{P}$) with $\delta =0$.
Multiplying in ($\mathcal{P}$) by $u_{i}$, integrating over $\Omega $,
applying Young inequality with $\alpha _{1},\beta _{2}>-1$, we get%
\begin{equation}
\begin{array}{l}
\int_{\Omega }\left\vert \nabla u\right\vert ^{p}dx=\lambda \int_{\Omega
}u^{\alpha _{1}+1}v^{\beta _{1}}dx\leq \lambda \int_{\Omega }(\frac{\alpha
_{1}+1}{p}u^{p}+\frac{p-1-\alpha _{1}}{p}v^{\frac{\beta _{1}p}{p-1-\alpha
_{1}}})\text{ }dx%
\end{array}
\label{55}
\end{equation}%
and%
\begin{equation}
\begin{array}{l}
\int_{\Omega }\left\vert \nabla v\right\vert ^{q}dx=\lambda \int_{\Omega
}u^{\alpha _{2}}v^{\beta _{2}+1}dx\leq \lambda \int_{\Omega }(\frac{%
q-1-\beta _{2}}{q}u^{\frac{\alpha _{2}q}{q-1-\beta _{2}}}+\frac{\beta _{2}+1%
}{q}v^{q})\text{ }dx.%
\end{array}
\label{56}
\end{equation}%
Adding (\ref{55}) with (\ref{56}), according to (\ref{c2}), this is
equivalent to%
\begin{equation}
\begin{array}{l}
\left\Vert \nabla u_{1}\right\Vert _{p_{1}}^{p_{1}}+\left\Vert \nabla
u_{2}\right\Vert _{p_{2}}^{p_{2}}\leq \lambda \lbrack (\frac{\alpha _{1}+1}{%
p_{1}}+\frac{p_{2}-1-\beta _{2}}{p_{2}})\left\Vert u_{1}\right\Vert
_{p_{1}}^{p_{1}}+(\frac{\beta _{2}+1}{p_{2}}+\frac{p_{1}-1-\alpha _{1}}{p_{1}%
})\left\Vert v\right\Vert _{p_{2}}^{p_{2}}].%
\end{array}
\label{60}
\end{equation}%
Since $\Theta =0$, observe from (\ref{c2}) that%
\begin{equation}
\left\{ 
\begin{array}{l}
\frac{\alpha _{1}+1}{p_{1}}+\frac{p_{2}-1-\beta _{2}}{p_{2}}=\frac{\alpha
_{1}+\alpha _{2}+1}{p_{1}} \\ 
\frac{\beta _{2}+1}{p_{2}}+\frac{p_{1}-1-\alpha _{1}}{p_{1}}=\frac{\beta
_{1}+\beta _{2}+1}{p_{2}}.%
\end{array}%
\right.  \label{61}
\end{equation}%
Then gathering (\ref{62}), (\ref{60}) and (\ref{61}) together yields%
\begin{equation*}
\begin{array}{c}
(\lambda _{1,p_{1}}-\frac{\alpha _{1}+\alpha _{2}+1}{p_{1}}\lambda
)\left\Vert u_{1}\right\Vert _{p_{1}}^{p_{1}}+(\lambda _{1,p_{2}}-\frac{%
\beta _{1}+\beta _{2}+1}{p_{2}}\lambda )\left\Vert u_{2}\right\Vert
_{p_{2}}^{p_{2}}\leq 0%
\end{array}%
\end{equation*}%
which is a contradiction for%
\begin{equation*}
\lambda <\lambda _{\ast }=\min \{\frac{p_{1}}{\alpha _{1}+\alpha _{2}+1}%
\lambda _{1,p_{1}},\frac{p_{2}}{\beta _{1}+\beta _{2}+1}\lambda _{1,p_{2}}\}.
\end{equation*}
Thus, problem ($\mathcal{P}$) has no solution for $\lambda <\lambda _{\ast
}, $ which completes the proof.

\end{document}